\documentclass[10pt]{article}
\title{\textbf{A Note on Semi-linear Wave Equations}}
\usepackage{amsmath,bbold}
\usepackage{geometry}
\usepackage{framed}
\usepackage{showkeys}
\usepackage{amsfonts,amssymb,amsmath,verbatim}
\usepackage{mathrsfs}
\usepackage{amssymb}
\usepackage{amsmath}
\usepackage{amsfonts}
\usepackage{mathrsfs}
\usepackage{amsfonts, amssymb, amsmath}
\usepackage{a4,graphicx}
\usepackage{graphicx}
\usepackage{slashed}
\usepackage{amsthm}
\usepackage{cite}
\newcommand{\leftexp}[2]{{\vphantom{#2}}^{#1}{#2}}
\author{Shuang Miao}
\begin{document}
\maketitle
\begin{abstract}
Inspired by the work of Wang and Yu \cite{WY1} on wave maps, we show that for all positive numbers $T_{0}>0$ and 
$E_{0}>0$, a large kind of semi-linear wave 
equation on $\mathbb{R}\times\mathbb{R}^{3}$ has a solution whose life-span is $[0,T_{0}]$, and the energy of the 
initial Cauchy data is at least $E_{0}$.
\end{abstract}
 \section{Introduction}
We consider in $\mathbb{R}\times\mathbb{R}^{3}$ the equation:
\begin{align}
 \Box\phi=\pm|\phi|^{k-1}\phi
\end{align}
where
\begin{align*}
 \Box=-\partial_{tt}^{2}+\Delta_{x},
\end{align*}
and $k$ is a odd number satisfying $k\geq3$\footnote{In this case, small data generates global solution, see \cite{Joh}. 
And we shall see that to establish the main estimates, we only need the nonlinearity to be a smooth function of $\phi$ 
which vanishes at the origin and whose growth is at least $|\phi|^{2}$ when $|\phi|$ is large.}. 
The equation (1) has a conserved energy
\begin{align}
 E(\phi(t))=\frac{1}{2}\int_{\mathbb{R}^{3}}|\partial_{t}\phi(t,x)|^{2}+|\nabla_{x}\phi(t,x)|^{2}\pm\frac{2}{k+1}|
\phi(t,x)|^{k+1}dx
\end{align}
The equation (1) is called defocusing, if there is a plus sign in front of the nonlinearity, otherwise it is called 
focusing. In view of the Sobolev embedding
$H^{1}(\mathbb{R}^{3})\hookrightarrow L^{6}(\mathbb{R}^{3})$, one refers to the range $k<5$ as the energy-subcritical 
regime, to $k=5$ as the energy-critical regime
and to $k>5$ as the energy-supercritical regime.
So the super-critical wave equations ($k>5$) (both defocusing and focusing case) are included in our present note.

    The study of the Cauchy problem for (1) has a long history. 
For the defocusing case, Rauch \cite{Rau} showed global existence for arbitrary 
smooth data for subcritical equations and for small energy smooth data in the critical case. Struwe \cite{Str1} 
obtained global existence for large but radially
symmetric data in the critical case, then Grillakis \cite{Gri} removed the radially symmetric condition on the data. 
After that Shatah and Struwe 
\cite{SS} studied the energy-class solutions. See also \cite{BG}, \cite{BS} and \cite{Tao1} for further results. 
The study of focusing case is initiated by Krieger 
and Schlag \cite{KS} as well as Kenig and Merle \cite{KM1}. Up to now, only a few results are known about the 
energy-supercritical case. See \cite{KM2}, 
\cite{KV} and \cite{Tao2}.

    Inspired by the recent work \cite{WY1}, we study long time solutions of semi-linear wave equations by a different 
approach. Following the ``short-pulse'' method, 
which was first introduced by Christodoulou \cite{Chr}, and extended by Klainerman and Rodnianski \cite{KR}, we 
establish the following long-time existence result for equation (1):

$\textbf{Main Theorem}$ \emph{For any $T_{0}>0$ and $E_{0}>0$, there exist $(\phi_{0},\phi_{1})\in C^{\infty}(\mathbb{R}^{3})
\times C^{\infty}(\mathbb{R}^{3})$ such that
the Cauchy problem for (1) with initial data $(\phi,\phi_{t})|_{t=0}=(\phi_{0},\phi_{1})$ has a unique solution $\phi\in 
C^{\infty}([0,T_{0}]\times\mathbb{R}^{3})$ with energy of at least $E_{0}$}.\\

     In \cite{WY1}, Wang and Yu constructed a solution for 2+1 wave maps with $\mathbb{S}^{2}$ as its target, by using 
a bootstrap argument. Since the characteristic initial data (so called $short pulse$) is chosen to be highly-oscillating,
 they can close the bootstrap. Then the solution will automatically have a life-span 
$[0,T_{0}]$, where $T_{0}$ is an arbitrary positive number given priorly. If the initial data is chosen properly, then 
the initial energy will be at least $E_{0}$, where $E_{0}$ is also an arbitrary positive number given priorly. The 
crucial point in their work is that the nonlinearity of wave maps into $\mathbb{S}^{2}$ is a 
``null form'', this means that the nonlinearity is ``not too bad'', so that they can absorb the term original from the 
nonlinearity in the a priori estimates if the characteristic initial data oscillates heavily enough. Wang and Yu also 
studied the 3+1 nonlinear wave equation with a ``null form'' by choosing the initial data at past null infinity, so 
they can even obtain global existence for large energy data, see \cite{WY2}. 

In the current work, there are no ``null form'' in the nonlinearity, but instead, the nonlinearity depends only on 
the solution itself, and we only commute
the ``bad'' vectorfield once with the operator $\Box$ when we do the bootstrap argument, so the nonlinearity will not 
cause trouble. Moreover, since the nonlinearity involves only the solution itself, we need one derivative less to close 
the bootstrap than the work of Wang and Yu. Our method for 3+1 semi-linear wave equations is also valid in the 2+1 case,
and we shall talk about this briefly at the end of the paper.

\section{Preliminaries}
\subsection{Basic Geometric Construction} 
\begin{figure} \label{cone}
 \centering
\includegraphics[width=0.9\textwidth,trim=0cm 15cm 0cm 8cm]{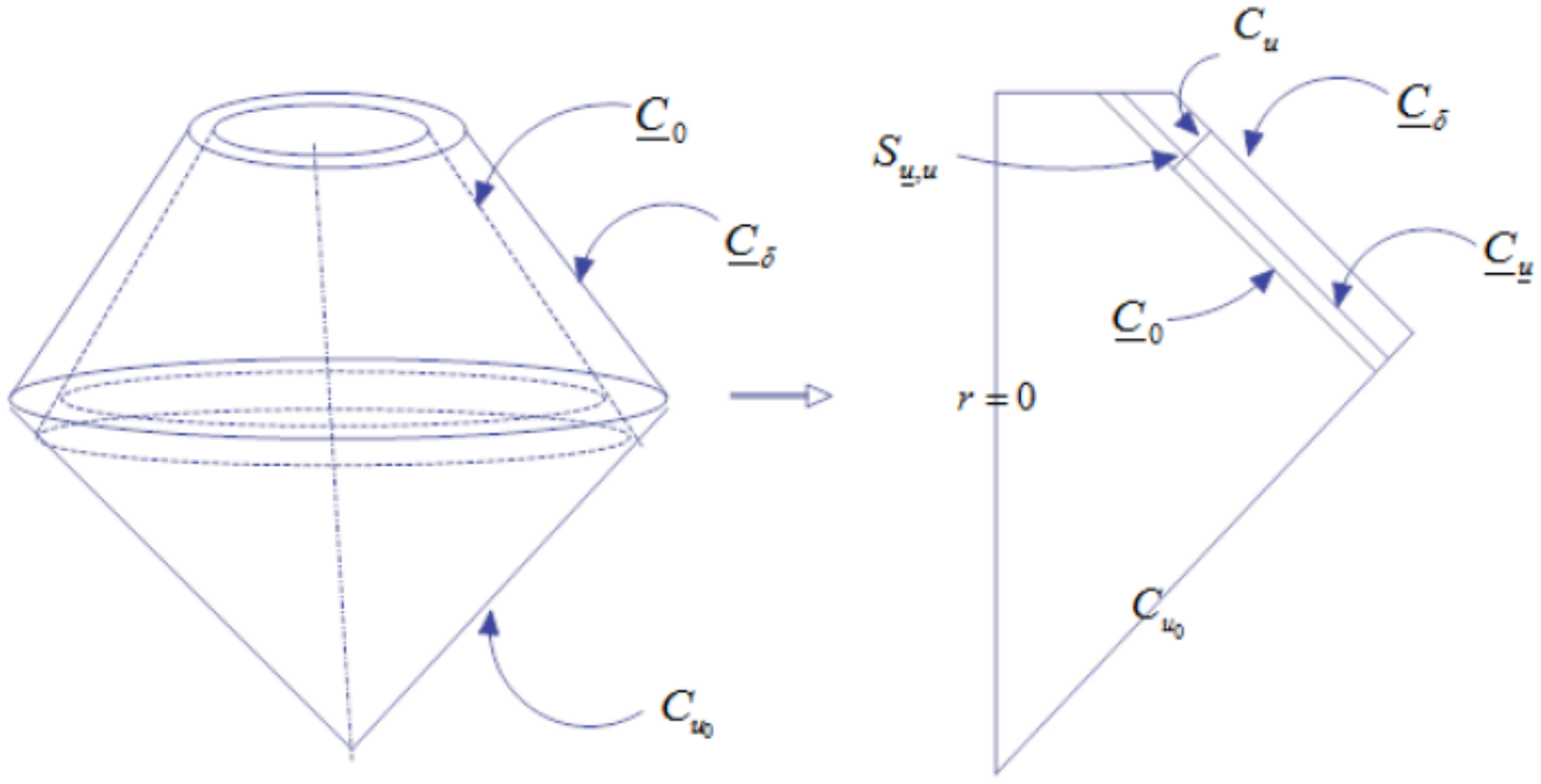} 
\caption{null cones (This figure is borrowed from \cite{WY1}.)} 
\end{figure}

This part is quite similar to \cite{WY1}, the
only difference is that  
we are in the 3-space dimensional Minkowski spacetime $\mathbb{R}\times\mathbb{R}^{3}$. 
We shall use the same notations as in \cite{WY1}.
 We have the optical functions:
\begin{align*}
 u=\frac{1}{2}(t-r)\quad \underline{u}=\frac{1}{2}(t+r)
\end{align*}
null vectorfields:
\begin{align*}
 L=\partial_{t}+\partial_{r}\quad\underline{L}=\partial_{t}-\partial_{r}
\end{align*}
as well as the rotation fields:
\begin{align*}
 \Omega_{ij}=x_{i}\partial_{j}-x_{j}\partial_{i}, 1\leq i,j\leq 3
\end{align*}
and also the relation between the operators $\Omega$ and $\slashed{\nabla}$:
\begin{align}
 |\Omega^{k}\phi|=|r|^{k}|\slashed{\nabla}^{k}\phi|
\end{align}
In section 3 and section 4, the parameter $u$ will be confined in the interval $[u_{0},-1]$ where $u_{0}\sim -T_{0}$ 
(then we can see that the life-span will be automatically $\sim T_{0}$). The parameter $\underline{u}$ is confined in 
$[0,\delta]$ where $\delta$ is a small parameter which will be determined later.
As in \cite{WY1}, the corresponding cones are pictured as follows (See \ref{cone}).
When we derive estimates in section 3 and 4, $\underline{u}\in[0,\delta]$ where $\delta$ will be sufficiently small. 
Since $T_{0}$ and $u_{0}$ are fixed numbers,
in the region where $(\underline{u},u]\in[0,\delta]\times[u_{0},-1]$, the parameter $r\sim 1$. In particular, we have
\begin{align*}
 |\Omega^{k}\phi|\sim |\slashed{\nabla}^{k}\phi|
\end{align*}
\subsection{Energy Identity}
Let $f$ be a solution for the following non-homogenous wave equation on $\mathbb{R}^{3+1}$:
\begin{align}
 \Box f=\Phi
\end{align}
The energy momentum tensor associated to $f$ is
\begin{align}
 T_{\mu\nu}[f]=\partial_{\mu}f\partial_{\nu}f-\frac{1}{2}g_{\mu\nu}|\nabla f|^{2}
\end{align}
Obviously, it is symmetric and satisfies the following identity:
\begin{align}
 \nabla^{\mu}T_{\mu\nu}[f]=\Phi\cdot\nabla_{\nu}f
\end{align}
Given a vectorfield $X$, which will be used as a $multiplier$ $vectorfield$, the associated energy currents are defined 
as follows
\begin{align*}
 P^{X}_{\alpha}[f]=T_{\alpha\mu}[f]X^{\mu},\quad K^{X}[f]=\frac{1}{2}T^{\mu\nu}[f]\leftexp{(X)}\pi_{\mu\nu}
\end{align*}
where the deformation tensor $\leftexp{(X)}{\pi}_{\mu\nu}$ is defined by
\begin{align}
 \leftexp{(X)}{\pi}_{\mu\nu}=\mathcal{L}_{X}g_{\mu\nu}=\nabla_{\mu}X_{\nu}
 +\nabla_{\nu}X_{\mu}
\end{align}
By (6), we easily obtain:
\begin{align}
 \nabla^{\mu}P^{X}_{\mu}[f]=K^{X}[f]+\Phi\cdot Xf
\end{align}

We can express $T$ in terms of null frames $\{e_{1}=r^{-1}\partial_{\theta}, e_{2}=(r\sin\theta)^{-1}, 
e_{3}=\underline{L}, e_{4}=L\}$:
\begin{align*}
 T(L,L)=|Lf|^{2}\quad T(\underline{L},\underline{L})=|\underline{L}f|^{2}\quad T(L,\underline{L})=|\slashed{\nabla}f|^{2}
\end{align*}
where we express the Minkowski metric in polar coordinates:
\begin{align*}
 g=-dt^{2}+dr^{2}+r^{2}(d\theta^{2}+(\sin\theta)^{2}d\varphi^{2})
\end{align*}

We shall use $X=L$ and $\underline{L}$ as multiplier vectorfields, the corresponding deformation tensors and currents 
are:
\begin{align}
  \leftexp{(L)}{\pi}=\frac{1}{r}\slashed{g}\quad\leftexp{(\underline{L})}{\pi}=-\frac{1}{r}\slashed{g}\\\notag
\quad K^{L}=\frac{1}{2r}Lf\underline{L}f,\quad K^{\underline{L}}=-\frac{1}{2r}Lf\underline{L}f
\end{align}
where $\slashed{g}$ is the restriction of the Minkowski metric to the sphere $S_{\underline{u},u}$.

\begin{figure} \label{domain}
 \centering 
\includegraphics[width=0.5\textwidth,trim=5cm 10cm 5cm 13cm]{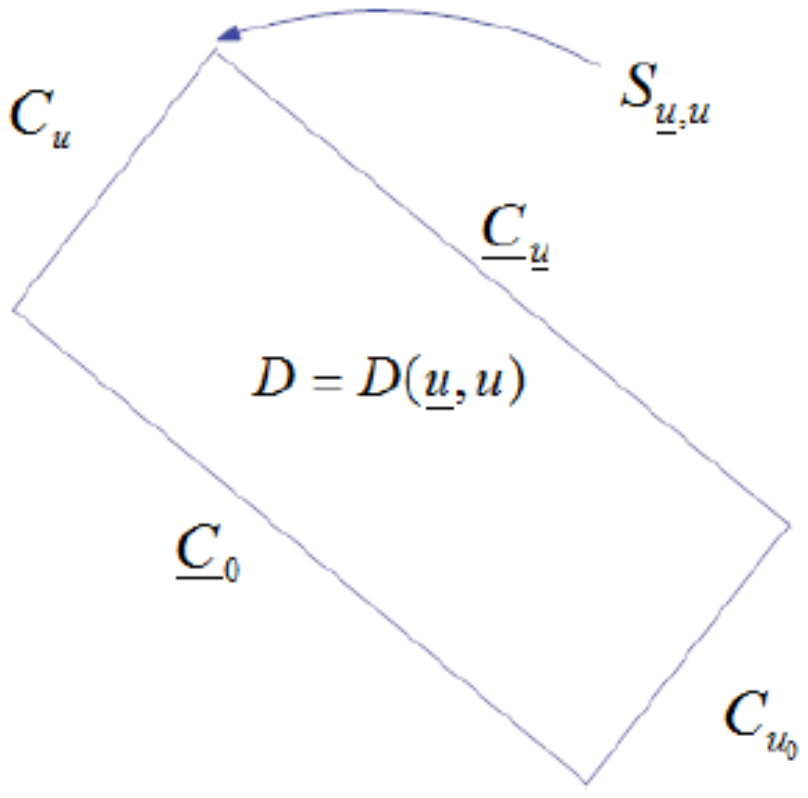} 
\caption{integral domain (This figure is borrowed from \cite{WY1}.)}
\end{figure}

We use $\mathcal{D}(u,\underline{u})$ to denote the space-time slab enclosed by the hypersurfaces $C_{u_{0}}$, 
$\underline{C}_{0}$, $C_{u}$ and $\underline{C}_{\underline{u}}$ as 
pictured above (see \ref{domain}). We integrate (8) on $\mathcal{D}(u,\underline{u})$ 
to obtain:
\begin{align*}
 \int_{C_{u}}T[f](X,L)+\int_{\underline{C}_{\underline{u}}}T[f](X,\underline{L})
=\int_{C_{u_{0}}}T[f](X,L)+\int_{\underline{C}_{0}}T[f](X,\underline{L})\\\notag
+\iint_{\mathcal{D}(u,\underline{u})}K^{X}[f]+\Phi\cdot Xf
\end{align*}
where $L$ and $\underline{L}$ are corresponding normals of the null hypersurfaces $C_{u}$ and 
$\underline{C}_{\underline{u}}$.

In applications, the data on $\underline{C}_{0}$ is always vanishing, thus, we have the following formula:
\begin{align}
 \int_{C_{u}}T[f](X,L)+\int_{\underline{C}_{\underline{u}}}T[f](X,\underline{L})
=\int_{C_{u_{0}}}T[f](X,L)+\iint_{\mathcal{D}(u,\underline{u})}K^{X}[f]+\Phi\cdot Xf
\end{align}

\subsection{Gronwall and Sobolev Inequalities}
We need the following Sobolev inequalities which can be derived from $isoperimetric$ $inequalty$. The proof can be 
found in \cite{Chr}.

$\textbf{Lemma 2.1}$ \emph{Let $(S,\slashed{g})$ be a compact 2-dimensional Riemannian manifold and $\phi$ a smooth function 
on $S$, which is square-integrable
and with square-integrable first derivatives. Then for $2<p<\infty$, $\phi\in L^{p}(S)$ and we have}:
\begin{align*}
 |S|^{-1/p}\|\phi\|_{L^{p}(S)}\leq C_{p}\sqrt{I^{\prime}(S)}\|\phi\|_{W^{2}_{1}(S)}
\end{align*}
\emph{Here $C_{p}$ is a numerical constant depending only on $p$, and $|S|$ is the area of the sphere $S_{\underline{u},u}$}.
\begin{align*}
 I^{\prime}(S)=\max\{I(S),1\}
\end{align*}
\emph{where $I(S)$ is the isoperimetric constant of $S$, and we define}:
\begin{align*}
 \|\phi\|_{W^{2}_{1}(S)}=\|\slashed{\nabla}\phi\|_{L^{2}(S)}+|S|^{-1/2}\|\phi\|_{L^{2}(S)}
\end{align*}

$\textbf{Lemma 2.2}$ \emph{Let $(S,\slashed{g})$ be a compact 2-dimensional Riemannian manifold and $\phi$ a smooth function 
on $S$, which belongs to $L^{p}(S)$ and with first derivatives which also belong to $L^{p}(S)$, for some $p>2$. 
Then $\phi\in L^{\infty}(S)$ and we have}
\begin{align*}
 \sup|\phi|\leq C_{p}\sqrt{I^{\prime}(S)}|S|^{(1/2)-(1/p)}\|\phi\|_{W^{p}_{1}(S)}
\end{align*}
\emph{Here $C_{p}$ is a numerical constant depending only on $p$, and we define}:
\begin{align*}
 \|\phi\|_{W^{p}_{1}(S)}=\|\slashed{\nabla}\phi\|_{L^{p}(S)}+|S|^{-1/2}\|\phi\|_{L^{p}(S)}
\end{align*}

Also by $isoperimetric$ $inequality$, we can deduce:

$\textbf{Lemma 2.3}$ \emph{Let $\phi$ be a smooth function on $C_{u}$ vanishing on $S_{0,u}$ then with the same condition as 
Lemma 2.1, we have}:
\begin{align*}
 \int_{S_{\underline{u}, u}}|\phi|^{6}d\mu_{\slashed{g}}\leq C(\int_{S_{\underline{u},u}}|\phi|^{4}d\mu_{\slashed{g}})
(|u|^{-2}\int_{S_{\underline{u},u}}|\phi|^{2}d\mu_{\slashed{g}}+\int_{S_{\underline{u},u}}|\slashed{\nabla}\phi|^{2}
d\mu_{\slashed{g}})
\end{align*}
\emph{where $C$ is an absolute constant}.

Obviously, from Lemma 2.1 and Lemma 2.2, we obtain:
\begin{align}
\sup_{S}|\phi|\leq C_{p}[|S|^{-1/2}\|\phi\|_{L^{2}(S)}+\|\slashed{\nabla}\phi\|_{L^{2}(S)}
+|S|^{1/2}\|\slashed{\nabla}^{2}\phi\|_{L^{2}(S)}] 
\end{align}

With the same condition as Lemma 2.3, we have:

$\textbf{Lemma 2.4}$ \emph{Let $\phi$ be a smooth $C_{u}$-function vanishing on $\underline{C}_{0}$, we have}:
\begin{align*}
 \sup_{\underline{u}}(|u|^{1/2}\|\phi\|_{L^{4}(S_{\underline{u},u})})\leq C_{p}\|L\phi\|^{1/2}_{L^{2}(C_{u})}
[\|\phi\|_{L^{2}(C_{u})}+|u|\|\slashed{\nabla}\phi\|_{L^{2}(C_{u})}]^{1/2}
\end{align*}
\emph{and by Gronwall's inequality, we have}:
\begin{align*}
 \|\phi\|_{L^{2}(S_{\underline{u},u})}\lesssim\|L\phi\|^{1/2}_{L^{2}(C_{u})}\|\phi\|^{1/2}_{L^{2}(C_{u})}
\end{align*}

Also from Lemma 2.3, we have:

$\textbf{Lemma 2.5}$ \emph{Let $\phi$ be a smooth function on $\underline{C}_{\underline{u}}$, the following estimates hold}:
\begin{align*}
 \sup_{u}(|u|^{1/2}\|\phi\|_{L^{4}(S_{\underline{u},u})})\leq C_{p}
\{|u_{0}|^{1/2}\|\phi\|_{L^{4}(S_{\underline{u},u_{0}})}+\||u|^{1/2}\underline{L}\phi\|^{1/2}
_{L^{2}(\underline{C}_{\underline{u}})}\cdot\\\
[\||u|^{-1/2}\phi\|_{L^{2}(\underline{C}_{\underline{u}})}^{2}+\||u|^{1/2}\slashed{\nabla}\phi\|^{2}
_{L^{2}(\underline{C}_{\underline{u}})}]^{1/4}\}
\end{align*}
\emph{also}:
\begin{align*}
 \|\phi\|_{L^{2}(S_{\underline{u},u})}\lesssim\|\phi\|_{L^{2}(S_{\underline{u},u_{0}})}
+\|\underline{L}\phi\|^{1/2}_{L^{2}(\underline{C}_{\underline{u}})}\|\phi\|^{1/2}_{L^{2}(\underline{C}_{\underline{u}})}
\end{align*}

From the above lemmas, we can easily obtain the following $L^{\infty}$ estimates:
\begin{align*}
 \|\phi\|_{L^{\infty}(S_{\underline{u},u})}\lesssim |u|^{-1/2}\|L\phi\|_{L^{2}(C_{u})}^{1/2}\|\phi\|
^{1/2}_{L^{2}(C_{u})}\\\notag
+\|L\slashed{\nabla}\phi\|_{L^{2}(C_{u})}^{1/2}\|\slashed{\nabla}\phi\|_{L^{2}(C_{u})}^{1/2}
+|u|^{1/2}\|L\slashed{\nabla}^{2}\phi\|_{L^{2}(C_{u})}^{1/2}\|\slashed{\nabla}^{2}\phi\|_{L^{2}(C_{u})}^{1/2}
\end{align*}
and also:
\begin{align*}
 \|\phi\|_{L^{\infty}(S_{\underline{u},u})}\lesssim |u|^{-1/2}\|\phi\|_{L^{2}(S_{\underline{u},u_{0}})}
+\|\slashed{\nabla}\phi\|_{L^{2}(S_{\underline{u},u_{0}})}
+|u|^{1/2}\|\slashed{\nabla}^{2}\phi\|_{L^{2}(S_{\underline{u},u_{0}})}\\\notag
+|u|^{-1/2}\|\underline{L}\phi\|_{L^{2}(\underline{C}_{\underline{u}})}^{1/2}\|\phi\|^{1/2}_{L^{2}(
\underline{C}_{\underline{u}})}
+\|\underline{L}\slashed{\nabla}\phi\|_{L^{2}(\underline{C}_{\underline{u}})}^{1/2}\|\slashed{\nabla}\phi\|
_{L^{2}(\underline{C}_{\underline{u}})}^{1/2}
+|u|^{1/2}\|\underline{L}\slashed{\nabla}^{2}\phi\|_{\underline{L}^{2}(\underline{C}_{\underline{u}})}^{1/2}
\|\slashed{\nabla}^{2}\phi\|_{L^{2}(\underline{C}_{\underline{u}})}^{1/2}
\end{align*}
In all of the above lemmas, we always assume that $\phi$ vanishes on $\underline{C}_{0}$.

Actually, in the following, we only use another version of the above Sobolev inequality, with 
$\slashed{\nabla}$ substituted by $\Omega$. In this case, 
the weight of $|u|$ will change, but it doesn't matter, because in our case, $|u|$ is more or less like a constant.\\

We also need the standard Gronwall's inequality:

$\textbf{Lemma 2.6}$ \emph{Let $f(t)$ be a non-negative function defined on an interval $I$ with initial point $t_{0}$. 
If $f$ satisfies}:
\begin{align*}
 \frac{d}{dt}f\leq a\cdot f+b
\end{align*}
\emph{where two non-negative functions $a, b\in L^{1}(I)$, then for all $t\in I$, we have}:
\begin{align}
 f(t)\leq e^{A(t)}(f(t_{0})+\int_{t_{0}}^{t}e^{-A(\tau)}b(\tau)d\tau)
\end{align}
\emph{where $A(t)=\int_{t_{0}}^{t}a(\tau)d\tau$}.

\subsection{Outline of the Proof}
We will follow the main steps of \cite{WY1}. The Cauchy data will be finally given on $t=u_{0}+\delta$ and the 
solution will exist at least for 
$t\in[u_{0}+\delta, -1]$. This can be shown in the above picture (see \ref{initialdata}):
\begin{figure}\label{initialdata}
 \centering 
\includegraphics[width=0.4\textwidth,trim=7cm 10cm 5cm 13cm]{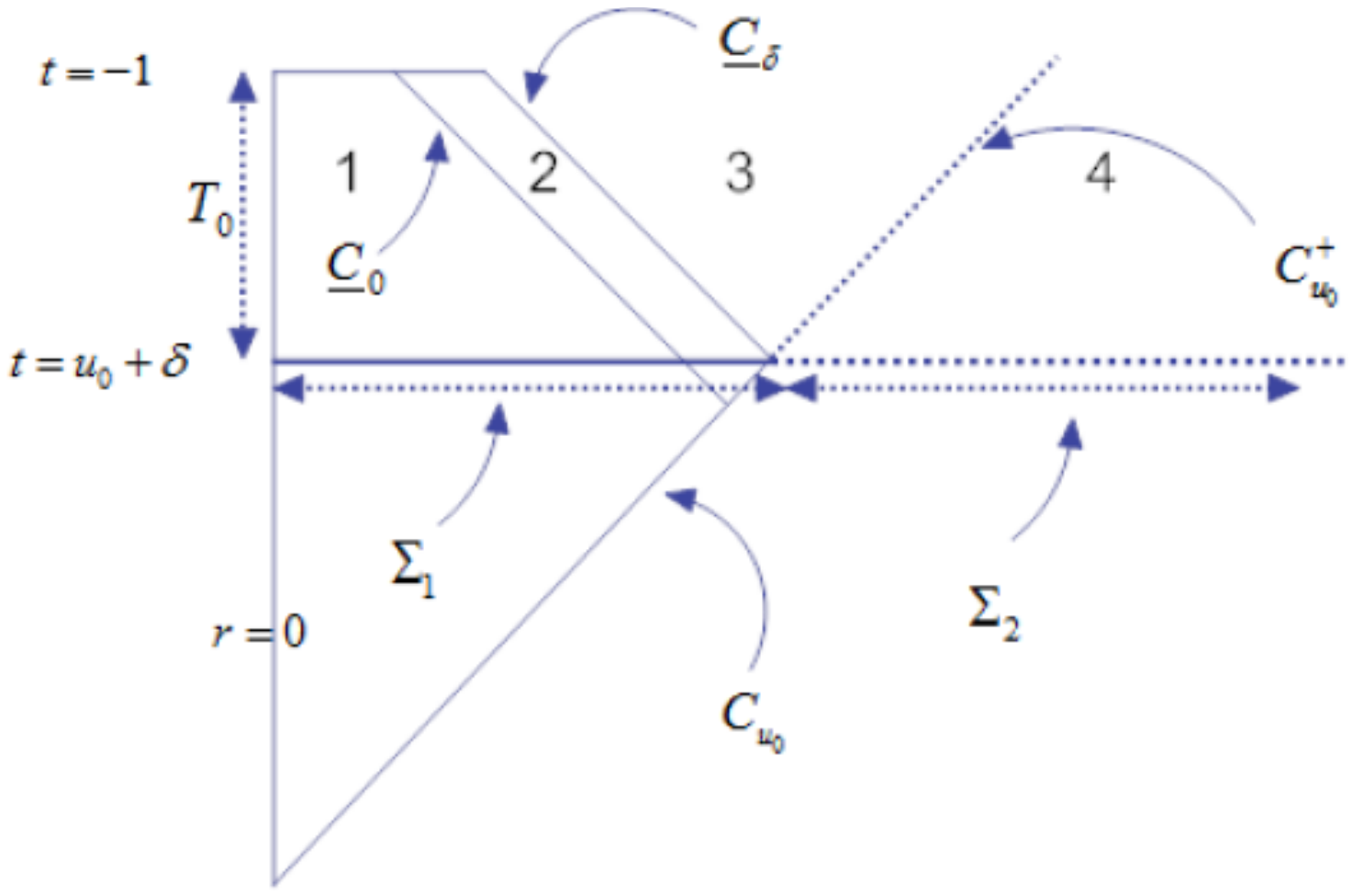}
\caption{construction of initial data (This figure is borrowed from \cite{WY1}.)}
\end{figure}
 First, we give initial data on the null hypersurface $C_{u_{0}}$ where $u_{0}\leq\underline{u}\leq\delta$. 
When $u_{0}\leq\underline{u}\leq 0$, the data is 
trivial, therefore the solution in Region 1 is zero. When $0\leq\underline{u}\leq\delta$, the data will be 
chosen as follows:
\begin{align*}
 \phi(\underline{u},u_{0},\theta)=\delta^{1/2}\psi_{0}(\frac{\underline{u}}{\delta},\theta)
\end{align*}
where the energy of $\psi_{0}$ is larger than $E_{0}$. We then show that we can construct a solution in Region 2. 
Consequently, we take the restriction of the solution 
constructed to the surface $\Sigma_{1}\subset\{t=u_{0}+\delta\}$ as the first part of the Cauchy data.\\

Second, we extend the Cauchy data on $\Sigma_{1}$ to $\Sigma_{2}\subset\{t=u_{0}+\delta\}$ such that the energy 
is small. By small data theory, we can construct a solution in 
Region 4.\\

Third, From previous two steps, we can show that the restriction of the solution already constructed to 
$\underline{C}_{\delta}$ and $C^{+}_{u_{0}}$ (where $
\underline{u}\geq\delta$) are small. We use them as initial data and we can solve this small data problem to 
construct solution on Region 3. We finally combines the solutions in Region 1,2,3 and 4 to finish the construction.

\section{Characteristic Initial Data}
    First, we require that the data $\phi(\underline{u},u_{0},\theta)$ to satisfy
\begin{align*}
 \phi(\underline{u},u_{0},\theta)=0\quad\textrm{for all}\quad \underline{u}\leq 0
\end{align*}
Therefore, according to the Huygens principle, the solution $\phi$ of (1) satisfies
\begin{align*}
 \phi\equiv0\quad\textrm{in Region 1}\quad=\{\underline{u}(x)\leq0, u_{0}\leq u(x)\leq 0\}
\end{align*}
Secondly, we choose
\begin{align}
 \phi(\underline{u},u_{0},\theta)=\delta^{1/2}\psi_{0}(\frac{\underline{u}}{\delta},\theta)
\end{align}
where $\psi_{0}$ is a smooth function supported in $(0,1)$ with respect to its first variable.

The data given in the above form is called a $short$ $pulse$, a name invented by Christodoulou in \cite{Chr}.\\

In order to derive the energy estimates, as in \cite{WY1}, we need the following commutators:
\begin{align}
 [L,\Omega]=0\quad [\underline{L},\Omega]=0\quad [\Box,\Omega]=0\\\notag
[\Box,L]=\frac{1}{r^{2}}(L-\underline{L})+\frac{2}{r^{3}}\slashed{\Delta}
\quad [\Box,\underline{L}]=\frac{1}{r^{2}}(\underline{L}-L)-\frac{2}{r^{3}}\slashed{\Delta}
\end{align}
Here the operator $\slashed{\Delta}$ is the Laplacian on standard sphere $\mathbb{S}^{2}$.

On the initial hypersurface $C_{u_{0}}$, we have the following bounds on data:
\begin{align*}
 \|L\phi\|_{L^{\infty}(C_{u_{0}})}\lesssim\delta^{-1/2}\quad \|\Omega\phi\|_{L^{\infty}(C_{u_{0}})}\lesssim\delta^{1/2}
\end{align*}
and for higher order derivatives, we have:
\begin{align*}
 \|L\Omega^{k}\phi\|_{L^{\infty}(C_{u_{0}})}\lesssim_{k}\delta^{-1/2}\\
\|\Omega^{k+1}\phi\|_{L^{\infty}(C_{u_{0}})}\lesssim_{k}\delta^{1/2}\\
\|L^{2}\Omega^{k-1}\phi\|_{L^{\infty}(C_{u_{0}})}\lesssim_{k}\delta^{-3/2}
\end{align*}
We also need a bound for $\underline{L}$-derivatives. To do this, we write the equation in null frames:
\begin{align}
 -L\underline{L}\phi+\frac{1}{r^{2}}\slashed{\Delta}\phi+\frac{1}{r}(L\phi-\underline{L}\phi)=\pm|\phi|^{k-1}\phi
\end{align}
We can write the above as a propagation equation for $\underline{L}\phi$ along $C_{u_{0}}$:
\begin{align*}
 L(\underline{L}\phi)=a\cdot\underline{L}\phi+b
\end{align*}
where 
\begin{align*}
 a=-\frac{1}{r}\quad b=\frac{1}{r}L\phi+\frac{1}{r^{2}}\slashed{\Delta}\phi\mp|\phi|^{k-1}\phi
\end{align*}
Obviously,
\begin{align*}
 \|a\|_{L^{\infty}(C_{u_{0}})}\lesssim 1,\quad \|b\|_{L^{\infty}(C_{u_{0}})}\lesssim\delta^{-1/2}
\end{align*}
Then by Gronwall's inequality, we easily obtain:
\begin{align}
 \|\underline{L}\phi\|_{L^{\infty}(C_{u_{0}})}\lesssim\delta^{1/2}
\end{align}
Similarly, by using the commutator, we obtain
\begin{align}
 \|\underline{L}\Omega\phi\|_{L^{\infty}(C_{u_{0}})}\lesssim\delta^{1/2}\quad\|\underline{L}\Omega^{2}\phi\|
_{L^{\infty}(C_{u_{0}})}\lesssim\delta^{1/2}
\end{align}
To obtain a long time existence theorem for (1), we have to derive estimates on $\phi$ as well as its derivatives. 
it's very natual that these estimates should 
be compatible with the bounds for $\phi$ on initial hypersurface. However, as stated in \cite{WY1}, a $relaxed$ 
estimate, which is easier to derive, is enough. That is, we just need the following bounds on $\Omega^{k}\phi$:
\begin{align}
 \|\Omega^{k+1}\phi\|_{L^{\infty}(C_{u_{0}})}\lesssim_{k}1
\end{align}
Summarizing, we have the following bounds on initial data:
\begin{align}
 \|L\Omega^{k}\phi\|_{L^{\infty}(C_{u_{0}})}\lesssim\delta^{-1/2}\\\notag
\|\Omega^{k+1}\phi\|_{L^{\infty}(C_{u_{0}})}\lesssim 1\\\notag
\|\underline{L}\Omega^{k}\phi\|_{L^{\infty}(C_{u_{0}})}\lesssim\delta^{1/2}
\end{align}
for $k=0,1,2$, and 
\begin{align}
 \|L^{2}\phi\|_{L^{\infty}(C_{u_{0}})}\lesssim\delta^{-3/2}\\\notag
\|L^{2}\Omega\phi\|_{L^{\infty}(C_{u_{0}})}\lesssim\delta^{-3/2}
\end{align}
From these $L^{\infty}$ bounds, we obtain easily the $L^{2}$ bounds:
\begin{align}
 \|L\Omega^{k}\phi\|_{L^{2}(C_{u_{0}})}\lesssim 1\\\notag
\|\Omega^{k+1}\phi\|_{L^{2}(C_{u_{0}})}\lesssim\delta^{1/2}
\end{align}
for $k=0,1,2$, and
\begin{align}
 \|L^{2}\phi\|_{L^{2}(C_{u_{0}})}\lesssim\delta^{-1}\\\notag
\|L^{2}\Omega\phi\|_{L^{2}(C_{u_{0}})}\lesssim\delta^{-1}
\end{align}
We shall show that (21) and (22) will hold on all later outgoing null hypersurfaces $C_{u}$ where $-1>u>u_{0}$ 
provided the solution of (1) can be constructed up to $C_{u}$.

\section{A priori Estimates}
    We start by defining a family of energy norms. For this purpose, we slightly abuse the notations: we use $C_{u}$ to 
denote $C^{[0,\underline{u}]}_{u}$ and $\underline{C}_{\underline{u}}$ to denote $\underline{C}_{\underline{u}}
^{[u_{0},u]}$, by definition,
\begin{align*}
 C^{[0,\underline{u}]}_{u}=\{p\in C_{u}|0\leq\underline{u}(p)\leq\underline{u}\}\quad
\underline{C}^{[u_{0},u]}_{\underline{u}}=\{p\in\underline{C}
_{\underline{u}}|u_{0}\leq u(p)\leq u\}
\end{align*}
We define the following norms which are the same order as in \cite{WY1}, but remember, we are now in $\mathbb{R}^{3+1}$ other than
$\mathbb{R}^{2+1}$. So essentially, we use one less derivative than 
\cite{WY1}. 
\begin{align}
 E_{1}(u,\underline{u})=\|L\phi\|_{L^{2}(C_{u})}+\delta^{-\frac{1}{2}}\|\Omega\phi\|_{L^{2}(C_{u})},\\\notag
\underline{E}_{1}(u,\underline{u})=\|\Omega\phi\|_{L^{2}(\underline{C}_{\underline{u}})}+
\delta^{-\frac{1}{2}}\|\underline{L}\phi\|_{L^{2}(\underline{C}_{\underline{u}})},\\\notag
E_{2}(u,\underline{u})=\|L\Omega\phi\|_{L^{2}(C_{u})}+\delta^{-\frac{1}{2}}\|\Omega^{2}\phi\|_{L^{2}(C_{u})},\\\notag
\underline{E}_{2}(u,\underline{u})=\|\Omega^{2}\phi\|_{L^{2}(\underline{C}_{\underline{u}})}+\delta^{-\frac{1}{2}}
\|\underline{L}\Omega\phi\|
_{L^{2}(\underline{C}_{\underline{u}})},\\\notag
E_{3}(u,\underline{u})=\|L\Omega^{2}\phi\|_{L^{2}(C_{u})}+\delta^{-\frac{1}{2}}\|\Omega^{3}\phi\|_{L^{2}(C_{u})},
\\\notag
\underline{E}_{3}(u,\underline{u})=\|\Omega^{3}\phi\|_{L^{2}(\underline{C}_{\underline{u}})}+\delta^{-\frac{1}{2}}
\|\underline{L}\Omega^{2}\phi\|
_{L^{2}(\underline{C}_{\underline{u}})}
\end{align}
 We also need another family of norms which involves at least two null derivatives. They are defined as follows:
\begin{align}
 F_{2}(u,\underline{u})=\delta\|L^{2}\phi\|_{L^{2}(C_{u})},\\\notag
\underline{F}_{2}(u,\underline{u})=\|\underline{L}^{2}\phi\|_{L^{2}(\underline{C}_{\underline{u}})},\\\notag
F_{3}(u,\underline{u})=\delta\|L^{2}\Omega\phi\|_{L^{2}(C_{u})},\\\notag
\underline{F}_{3}(u,\underline{u})=\|\underline{L}^{2}\Omega\phi\|_{L^{2}(\underline{C}_{\underline{u}})}.
\end{align}

We shall prove as in \cite{WY1}:

$\textbf{Main A priori Estimates}$. If $\delta$ is sufficiently small, for all initial data of (1) and all 
$I\in\mathbb{R}_{>0}$ which satisfy
\begin{align}
 E_{1}(u_{0},\delta)+E_{2}(u_{0},\delta)+E_{3}(u_{0},\delta)+F_{2}(u_{0},\delta)+F_{3}(u_{0},\delta)\leq I,
\end{align}
there is a constant $C(I)$ depending only on $I$ (in particular, not on $\delta$), so that
\begin{align}
 \sum_{i=1}^{3}[E_{i}(u,\underline{u})+\underline{E}_{i}(u,\underline{u})]+\sum_{j=2}^{3}[F_{j}(u,\underline{u})
+\underline{F}_{j}(u,\underline{u})]\leq C(I),
\end{align}
for all $u\in[u_{0},u^{*}]$ and $\underline{u}\in[0,\underline{u}^{*}]$ where $u_{0}\leq u^{*}\leq -1$ and $0\leq 
\underline{u}^{*}\leq \delta$.

 We consider the set $\mathcal{A}\subset\mathcal{U}:
=\{(u,\underline{u}): u\in[u_{0},u^{*}], \underline{u}\in[0,\underline{u}^{*}]\}$, in which the following holds:
\begin{align}
M:=\sum_{i=1}^{3}[E_{i}(u,\underline{u})+\underline{E}_{i}(u,\underline{u})]+\sum_{j=2}^{3}[F_{j}(u,\underline{u})
+\underline{F}_{j}(u,\underline{u})]\leq AI,
\end{align}
 where $A>1$ is a sufficiently large constant depending only on initial data. Obviously, $\mathcal{A}$ is not empty. 
Here $\mathcal{U}$ is the set where the solution exists. We shall prove that actually $\mathcal{A}=\mathcal{U}$.

\subsection{Preliminary Estimates} Under the bootstrap assumption (27), we first derive $L^{\infty}$ for one derivatives 
of $\phi$. We will also obtain the $L^{4}$ estimates for derivatives of $\phi$ up to the second order. 

     We start with $L\phi$. According to Sobolev inequalities, we have (we shall omit the weights on $|u|$)
\begin{align*}
 \|L\phi\|_{L^{4}(S_{\underline{u},u})}\lesssim\|L^{2}\phi\|^{1/2}_{L^{2}(C_{u})}
(\|L\phi\|^{1/2}_{L^{2}(C_{u})}+\|L\Omega\phi\|^{1/2}_{L^{2}(C_{u})})\\
\lesssim(\delta^{-1}M)^{\frac{1}{2}}M^{\frac{1}{2}}
\end{align*}
Hence,
\begin{align}
 \|L\phi\|_{L^{4}(S_{\underline{u},u})}\lesssim\delta^{-1/2}M
\end{align}
Similarly, we have
\begin{align}
 \|L\Omega\phi\|_{L^{4}(S_{\underline{u},u})}\lesssim\delta^{-1/2}M
\end{align}

Now we consider $\Omega\phi$. According to Sobolev inequalities, we have
\begin{align*}
 \|\Omega\phi\|_{L^{4}(S_{\underline{u},u})}\lesssim\|L\Omega\phi\|_{L^{2}(C_{u})}^{1/2}
(\|\Omega\phi\|^{1/2}_{L^{2}(C_{u})}+\|\Omega^{2}\phi\|^{1/2}_{L^{2}(C_{u})})\\
\lesssim M^{1/2}((\delta^{1/2}M)^{1/2}+(\delta^{1/2}M)^{1/2})
\end{align*}
Thus,
\begin{align}
 \|\Omega\phi\|_{L^{4}(S_{\underline{u},u})}\lesssim\delta^{1/4}M
\end{align}
Similarly,
\begin{align}
 \|\Omega^{2}\phi\|_{L^{4}(S_{\underline{u},u})}\lesssim\delta^{1/4}M
\end{align}
 Finally, we turn to the estimates on $\underline{L}\phi$.
\begin{align*}
 \|\underline{L}\phi\|_{L^{4}(S_{\underline{u},u})}\lesssim\|\underline{L}\phi\|_{L^{4}(S_{\underline{u},u_{0}})}
+\|\underline{L}^{2}\phi\|_{L^{2}(\underline{C}_{\underline{u}})}^{1/2}(\|\underline{L}\phi\|_{L^{2}(\underline{C}
_{\underline{u}})}^{1/2}+\|\Omega\underline{L}\phi\|_{L^{2}(\underline{C}_{\underline{u}})}^{1/2})\\
\lesssim\delta^{1/2}+\delta^{1/4}M
\end{align*}
If $\delta$ is sufficiently small, we obtain
\begin{align}
 \|\underline{L}\phi\|_{L^{4}(S_{\underline{u},u})}\lesssim\delta^{1/4}M
\end{align}
Similarly, we also obtain
\begin{align}
 \|\underline{L}\Omega\phi\|_{L^{4}(S_{\underline{u},u})}\lesssim\delta^{1/4}M
\end{align}
Finally we need an $L^{\infty}$ bound for $\phi$. If we set 
\begin{align*}
 \|\phi\|_{L^{\infty}(S_{\underline{u},u})}\leq N
\end{align*}
then by the first $L^{\infty}$ Sobolev inequality, we have:
\begin{align*}
 N\lesssim\delta^{1/4}N^{1/2}M^{1/2}+\delta^{1/4}M
\end{align*}
So if $\delta$ is sufficiently small, we have:
\begin{align}
 \|\phi\|_{L^{\infty}(S_{\underline{u},u})}\lesssim\delta^{1/4}M
\end{align}
We summarize all the estimates in the following proposition.

$\textbf{Proposition 4.1}$ Under the bootstrap assumption (27), if $\delta$ is sufficiently small, we have
\begin{align*}
 \delta^{-1/4}\|\phi\|_{L^{\infty}(S_{\underline{u},u})}
+\delta^{1/2}\|L\Omega\phi\|_{L^{4}(S_{\underline{u},u})}+\delta^{-1/4}\|\Omega^{2}\phi\|_{L^{4}(S_{\underline{u},u})}
+\delta^{-1/4}\|\underline{L}\Omega\phi\|_{L^{4}(S_{\underline{u},u})}\\
+\delta^{1/2}\|L\phi\|_{L^{4}(S_{\underline{u},u})}+\delta^{-1/4}\|\Omega\phi\|_{L^{4}(S_{\underline{u},u})}
+\delta^{-1/4}\|\underline{L}\phi\|_{L^{4}(S_{\underline{u},u})}\lesssim M
\end{align*}

\subsection{Estimates on $E_{k}$ and $\underline{E}_{k}$} For simplicity, we shall assume that $k=2$, because if $k>2$, 
one just need to bound the extra power by $L^{\infty}$ norm.

We commute $\Omega^{i}$ (for $i=1,2$) with (1), we have \footnote{For simplicity, we omit the constant coefficients 
and the sign for nonlinearity.}
\begin{align*}
 \Box\Omega^{i}\phi=\phi\Omega^{k}\phi+\sum_{|p|>0,|q|>0,p+q=i}\Omega^{p}\phi\Omega^{q}\phi
\end{align*}
Now we use the basic energy identity for this equation where we take $f=\Omega^{i}\phi$ $(i=0,1,2)$ and $X=L$, then we 
have:
\begin{align}
 \int_{C_{u}}|L\Omega^{i}\phi|^{2}+\int_{\underline{C}_{\underline{u}}}|\slashed{\nabla}\Omega^{i}\phi|^{2}
=\int_{C_{u_{0}}}|L\Omega^{i}\phi|^{2}+\\\notag
\iint_{\mathcal{D}}(\phi\Omega^{i}\phi)L\Omega^{i}\phi+\sum_{|p|>0,|q|>0,p+q=i}\iint_{\mathcal{D}}
(\Omega^{p}\phi\Omega^{q}\phi)
L\Omega^{i}\phi+\iint_{\mathcal{D}}\frac{1}{2r}\underline{L}\Omega^{i}\phi L\Omega^{i}\phi\\\notag
=\int_{C_{u_{0}}}|L\Omega^{i}\phi|^{2}+S_{1}+S_{2}+S_{3}
\end{align}
where $S_{j}$ are defined in the obvious way. We also recall that:
\begin{align*}
 \|\Omega^{i}\phi\|_{L^{p}(S_{\underline{u},u})}\sim\|\slashed{\nabla}^{i}\phi\|_{L^{p}(S_{\underline{u},u})}
\end{align*}
By Proposition 4.1 and bootstrap assumption,
\begin{align*}
 S_{1}\lesssim\|\phi\|_{L^{\infty}(S_{\underline{u},u})}\int_{u_{0}}^{u}\|\Omega^{i}\phi\|_{L^{2}(C_{u^{\prime}})}
\|L\Omega^{i}\phi\|_{L^{2}(C_{u^{\prime}})}du^{\prime}
\lesssim(\delta^{1/4}M)(\delta^{1/2}M)M\lesssim\delta^{3/4}M^{3}
\end{align*}
also
\begin{align*}
 S_{2}\lesssim\int_{u_{0}}^{u}\|\Omega^{p}\phi\|_{L^{4}(C_{u^{\prime}})}\|\Omega^{q}\phi\|_{L^{4}(C_{u^{\prime}})}
\|L\Omega^{i}\phi\|_{L^{2}(C_{u^{\prime}})}du^{\prime}\lesssim\int_{u_{0}}^{u}(\delta^{1/2}M)(\delta^{1/2}M)Mdu^{\prime}
\lesssim\delta M^{3}
\end{align*}
and
\begin{align*}
S_{3}\lesssim\|L\Omega^{i}\phi\|_{L^{2}(\mathcal{D})}\|\underline{L}\Omega^{i}\phi\|_{L^{2}(\mathcal{D})}\\
\lesssim(\int_{u_{0}}^{u}\|L\Omega^{i}\phi\|^{2}_{L^{2}(C_{u^{\prime}})}du^{\prime})^{1/2}(\int_{0}^{\underline{u}}\|
\underline{L}\Omega^{i}\phi\|^{2}_{L^{2}(\underline{C}_{\underline{u}^{\prime}})}d\underline{u}^{\prime})^{1/2}
\lesssim M(\delta M)\lesssim\delta M^{2}
\end{align*}
Put all these in (35), we obtain:
\begin{align}
 \sum_{i=0}^{2}(\|L\Omega^{i}\varphi\|_{L^{2}(C_{u})}+\|\slashed{\nabla}\Omega^{i}\varphi\|_{L^{2}(\underline{C}_{\underline{u}})})\lesssim
I+\delta^{3/8}M^{3/2}
\end{align}
We still consider, for $i=0,1,2$,
\begin{align*}
 \Box\Omega^{i}\phi=\phi\Omega^{k}\phi+\sum_{|p|>0,|q|>0,p+q=i}\Omega^{p}\phi\Omega^{q}\phi
\end{align*}
But now we take $X=\underline{L}$ in the energy identity:
\begin{align}
 \int_{C_{u}}|\slashed{\nabla}\Omega^{i}\phi|^{2}+\int_{\underline{C}_{\underline{u}}}|\underline{L}\Omega^{i}\phi|^{2}
=\int_{C_{u_{0}}}|\slashed{\nabla}\Omega^{i}\phi|^{2}\\\notag
+\iint_{\mathcal{D}}\phi\Omega^{i}\phi \underline{L}\Omega^{i}\phi+\sum_{|p|>0,|q|>0,p+q=i}\iint_{\mathcal{D}}\Omega^{p}\phi\Omega^{q}\phi
(\underline{L}\Omega^{i}\phi)-\iint_{\mathcal{D}}\frac{1}{2r}(\underline{L}\Omega^{i}\phi)(L\Omega^{i}\phi)\\\notag
=\int_{C_{u_{0}}}|\slashed{\nabla}\Omega^{i}\phi|^{2}+T_{1}+T_{2}+T_{3}
\end{align}
where $T_{j}$ are defined in the obvious way.

As before, by Proposition 4.1 and bootstrap assumption, we have:
\begin{align*}
 T_{1}\lesssim\|\phi\|_{L^{\infty}(S_{\underline{u},u})}\int_{0}^{\underline{u}}\|\Omega^{i}\phi\|_{L^{2}(\underline{C}
_{\underline{u}^{\prime}})}
\|\underline{L}\Omega^{i}\phi\|_{L^{2}(\underline{C}_{\underline{u}^{\prime}})}d\underline{u}^{\prime}\lesssim
(\delta^{1/4}M)(M)(\delta^{1/2}M)\delta\lesssim\delta^{7/4}M^{3}
\end{align*}
also
\begin{align*}
 T_{2}\lesssim\int_{0}^{\underline{u}}\|\Omega^{p}\phi\|_{L^{4}(\underline{C}_{\underline{u}^{\prime}})}
\|\Omega^{q}\phi\|_{L^{4}(\underline{C}_{\underline{u}^{\prime}})}\|\underline{L}\Omega^{i}\phi\|_{L^{2}(\underline{C}
_{\underline{u}^{\prime}})}
d\underline{u}^{\prime}\lesssim(\delta^{1/4}M)(\delta^{1/4}M)(\delta^{1/2}M)\delta\lesssim\delta^{2}M^{3}
\end{align*}
Similar to $S_{3}$, we have:
\begin{align*}
 T_{3}\lesssim\delta M^{2}
\end{align*}
Put these in (37), we obtain:
\begin{align}
 \sum_{i=0}^{2}(\|\slashed{\nabla}\Omega^{i}\varphi\|_{L^{2}(C_{u})}+\|\underline{L}\Omega^{i}\varphi\|
_{L^{2}(\underline{C}_{\underline{u}})})
\lesssim I+\delta^{1/2}M^{3/2}
\end{align}
Combining (36) and (38) we obtain:
\begin{align}
 \sum_{i=1}^{3}(E_{i}(u,\underline{u})+\underline{E}_{i}(u,\underline{u}))\lesssim I+\delta^{3/8}M^{3/2}
\end{align}

\subsection{Estimates on $F_{2}(u,\underline{u})$ and $\underline{F}_{2}(u,\underline{u})$} We first consider the bound 
of $\|\underline{L}^{2}\varphi\|_{L^{2}(\underline{C}_{\underline{u}})}$. We commute $\underline{L}$ with (1), we obtain
\begin{align*}
 \Box\underline{L}\phi=\phi\underline{L}\phi+\frac{1}{r^{2}}(\underline{L}\phi-L\phi)-\frac{2}{r^{3}}\slashed{\Delta}\phi
\end{align*}
We use the basic energy identity where we take $f=\underline{L}\phi$ and $X=\underline{L}$, therefore,
\begin{align}
 \int_{C_{u}}|\slashed{\nabla}\underline{L}\phi|^{2}+\int_{\underline{C}_{\underline{u}}}|\underline{L}^{2}\phi|^{2}
=\int_{C_{u_{0}}}|\slashed{\nabla}\underline{L}\phi|^{2}+\iint_{\mathcal{D}}\phi\underline{L}\phi\underline{L}^{2}\phi\\
+\iint_{\mathcal{D}}\frac{1}{r^{2}}(\underline{L}\phi-L\phi)\underline{L}^{2}\phi+\iint_{\mathcal{D}}\frac{1}{r^{3}}
\slashed{\Delta}\phi\underline{L}^{2}\phi\\
=\int_{C_{u_{0}}}|\slashed{\nabla}\underline{L}\phi|^{2}+S_{1}+S_{2}+S_{3}
\end{align}
For $S_{1}$, we have, by Proposition 4.1 and bootstrap assumption:
\begin{align*}
 S_{1}\lesssim\|\phi\|_{L^{\infty}(S_{\underline{u},u})}\int_{0}^{\underline{u}}\|\underline{L}\phi\|
_{L^{2}(\underline{C}_{\underline{u}^{\prime}})}
\|\underline{L}^{2}\phi\|_{L^{2}(\underline{C}_{\underline{u}^{\prime}})}d\underline{u}^{\prime}
\lesssim(\delta^{1/4}M)(\delta^{1/2}M)(M)\delta\lesssim\delta^{5/4}M^{3}
\end{align*}
also,
\begin{align*}
 S_{2}\lesssim\int_{0}^{\underline{u}}\|\underline{L}\phi\|_{L^{2}(\underline{C}_{\underline{u}^{\prime}})}
\|\underline{L}^{2}\phi\|_{L^{2}(\underline{C}_{\underline{u}^{\prime}})}d\underline{u}^{\prime}+
(\int_{u_{0}}^{u}\|L\phi\|^{2}_{L^{2}(C_{u^{\prime}})}du^{\prime})^{1/2}(\int_{0}^{\underline{u}}
\|\underline{L}^{2}\phi\|^{2}_{L^{2}(\underline{C}
_{\underline{u}^{\prime}})}d\underline{u}^{\prime})^{1/2}\\
\lesssim\delta^{3/2}M^{2}+\delta^{1/2}M^{2}\lesssim\delta^{1/2}M^{2}
\end{align*}
The estimates for $S_{3}$ is similar,
\begin{align*}
 S_{3}\lesssim\int_{0}^{\underline{u}}\|\Omega^{2}\phi\|_{L^{2}(\underline{C}_{\underline{u}^{\prime}})}\|
\underline{L}^{2}\phi\|_{L^{2}(\underline{C}_{\underline{u}^{\prime}})}d\underline{u}^{\prime}\lesssim\delta M^{2}
\end{align*}
So we obtain:
\begin{align}
 \underline{F}_{2}(u,\underline{u})\lesssim\delta I+\delta^{1/4}M^{3/2}
\end{align}

Next, we consider the bound for $\|L^{2}\phi\|_{L^{2}(C_{u})}$, we commute $L$ with (1):
\begin{align*}
 \Box L\phi=\phi L\phi-\frac{1}{r^{2}}(\underline{L}\phi-L\phi)+\frac{2}{r^{3}}\slashed{\Delta}\phi
\end{align*}
We use the energy identity with $f=L\phi$ and $X=L$,
\begin{align*}
 \int_{C_{u}}|L^{2}\phi|^{2}+\int_{\underline{C}_{\underline{u}}}|\Omega L\phi|^{2}=\int_{C_{u_{0}}}|L^{2}\phi|^{2}
+\iint_{\mathcal{D}}\phi L\phi L^{2}\phi\\
+\iint_{\mathcal{D}}\frac{1}{r^{2}}(\underline{L}\phi-L\phi)L^{2}\phi
+\iint_{\mathcal{D}}\frac{1}{r^{3}}\slashed{\Delta}\phi L^{2}\phi\\
=\int_{C_{u_{0}}}|L^{2}\phi|^{2}+T_{1}+T_{2}+T_{3}
\end{align*}
As usual, by Proposition 4.1 and bootstrap assumption, we have:
\begin{align*}
 T_{1}\lesssim\|\phi\|_{L^{\infty}(S_{\underline{u},u})}\int_{u_{0}}^{u}\|L\phi\|_{L^{2}(C_{u^{\prime}})}\|L^{2}\phi\|
_{L^{2}(C_{u^{\prime}})}du^{\prime}
\lesssim(\delta^{1/4}M)M(\delta^{-1}M)\lesssim\delta^{-3/4}M^{3}
\end{align*}
also,
\begin{align*}
 T_{2}\lesssim(\int_{0}^{\underline{u}}\|\underline{L}\phi\|^{2}_{L^{2}(\underline{C}_{\underline{u}^{\prime}})}
d\underline{u}^{\prime})^{1/2}
(\int_{u_{0}}^{u}\|L^{2}\phi\|^{2}_{L^{2}(C_{u^{\prime}})}du^{\prime})^{1/2}+\int_{u_{0}}^{u}\|L\phi\|
_{L^{2}(C_{u^{\prime}})}\|L^{2}\phi\|_{L^{2}(C_{u^{\prime}})}du^{\prime}\\
\lesssim\delta^{-1}M^{2}
\end{align*}
Estimates for $T_{3}$ is similar,
\begin{align*}
 T_{3}\lesssim\int_{u_{0}}^{u}\|\Omega^{2}\phi\|_{L^{2}(C_{u^{\prime}})}\|L^{2}\phi\|_{L^{2}(C_{u^{\prime}})}du^{\prime}
\lesssim(\delta^{1/2}M)(\delta^{-1}M)\lesssim\delta^{-1/2}M^{2}
\end{align*}
So we have:
\begin{align*}
 \int_{C_{u}}|L^{2}\phi|^{2}\lesssim\delta^{-2}I^{2}+\delta^{-1}M^{3}
\end{align*}
thus,
\begin{align}
 F_{2}(u,\underline{u})=\delta\|L^{2}\phi\|_{L^{2}(C_{u})}\lesssim I+\delta^{1/2}M^{3/2}
\end{align}
Summarizing, we obtain:
\begin{align}
 F_{2}(u,\underline{u})\lesssim I+\delta^{1/2}M^{3/2}\\\notag
\underline{F}_{2}(u,\underline{u}))\lesssim\delta I+\delta^{1/4}M^{3/2}
\end{align}
\subsection{Estimates on $F_{3}(u,\underline{u})$ and $\underline{F}_{3}(u,\underline{u})$} 
We first estimate $F_{3}(\underline{u},u)$, we commute $L$ and 
$\Omega$ with (1) to derive
\begin{align*}
 \Box L\Omega\phi=\phi L\Omega\phi+\Omega\phi L\phi-\frac{1}{r^{2}}(\underline{L}\Omega\phi-L\Omega\phi)
+\frac{2}{r^{3}}\slashed{\Delta}\Omega\phi
\end{align*}
Applying the energy identity with $f=L\Omega\phi$, and $X=L$, we obtain:
\begin{align*}
 \int_{C_{u}}|L^{2}\Omega\phi|^{2}+\int_{\underline{C}_{\underline{u}}}|\slashed{\nabla}L\Omega\phi|^{2}=
\int_{C_{u_{0}}}|L^{2}\Omega\phi|^{2}+\iint_{\mathcal{D}}\phi L\Omega\phi L^{2}\Omega\phi+
\iint_{\mathcal{D}}\Omega\phi L\phi L^{2}\Omega\phi\\
+\iint_{\mathcal{D}}\frac{1}{r^{2}}(\underline{L}\Omega\phi-L\Omega\phi)L^{2}\Omega\phi+\iint_{\mathcal{D}}
\frac{1}{r^{3}}
\slashed{\Delta}\Omega\phi L^{2}\Omega\phi\\
=\int_{C_{u_{0}}}|L^{2}\Omega\phi|^{2}+S_{1}+S_{2}+S_{3}+S_{4}
\end{align*}
As before, by Proposition 4.1 and bootstrap assumption, we have:
\begin{align*}
 S_{1}\lesssim\|\phi\|_{L^{\infty}(S_{\underline{u},u})}\int_{u_{0}}^{u}\|L\Omega\phi\|_{L^{2}(C_{u^{\prime}})}
\|L^{2}\Omega\phi\|_{L^{2}(C_{u^{\prime}})}du^{\prime}\lesssim(\delta^{1/4}M)M(\delta^{-1}M)\lesssim\delta^{-3/4}M^{3}
\end{align*}
also
\begin{align*}
 S_{2}\lesssim\int_{u_{0}}^{u}\|\Omega\phi\|_{L^{4}(C_{u^{\prime}})}\|L\phi\|_{L^{4}(C_{u^{\prime}})}\|L^{2}\Omega\phi\|
_{L^{2}(C_{u^{\prime}})}du^{\prime}\lesssim(\delta^{1/2}M)(\delta^{1/2}M)(\delta^{-1}M)\lesssim M^{3}
\end{align*}
For $S_{3}$, we have:
\begin{align*}
 S_{3}\lesssim(\int_{0}^{\underline{u}}\|\underline{L}\Omega\phi\|^{2}_{L^{2}(C_{\underline{u}^{\prime}})}
d\underline{u}^{\prime})^{1/2}(\int_{u_{0}}^{u}\|L^{2}\Omega\phi\|^{2}_{L^{2}(C_{u^{\prime}})}du^{\prime})^{1/2}+
\int_{u_{0}}^{u}\|L\Omega\phi\|_{L^{2}(C_{u^{\prime}})}\|L^{2}\Omega\phi\|_{L^{2}(C_{u^{\prime}})}du^{\prime}\\
\lesssim(\delta M)(\delta^{-1}M)+M(\delta^{-1}M)\lesssim\delta^{-1}M^{2}
\end{align*}
Similarly,
\begin{align*}
 S_{4}\lesssim\int_{u_{0}}^{u}\|\slashed{\Delta}\Omega\phi\|_{L^{2}(C_{u^{\prime}})}\|L^{2}\Omega\phi\|
_{L^{2}(C_{u^{\prime}})}du^{\prime}\lesssim(\delta^{1/2}M)(\delta^{-1}M)\lesssim\delta^{-1/2}M^{2}
\end{align*}
So we obtain:
\begin{align*}
 \int_{C_{u}}|L^{2}\Omega\phi|^{2}\lesssim\delta^{-2}I^{2}+\delta^{-1}M^{3}
\end{align*}
i.e.
\begin{align}
 F_{3}(u,\underline{u})\lesssim I+\delta^{1/2}M^{3/2}
\end{align}

For $\underline{F}_{3}(\underline{u},u)$, we commute $\underline{L}$ and $\Omega$ with (1) to 
derive
\begin{align*}
 \Box\underline{L}\Omega\phi=\phi\underline{L}\Omega\phi+\underline{L}\phi\Omega\phi
+\frac{1}{r^{2}}(\underline{L}\Omega\phi-L\Omega\phi)-\frac{2}{r^{3}}\slashed{\Delta}\Omega\phi
\end{align*}
We use the energy identity with $f=\underline{L}\Omega\phi$ and $X=\underline{L}$,
\begin{align*}
 \int_{C_{u}}|\slashed{\nabla}\underline{L}\Omega\phi|^{2}+\int_{\underline{C}_{\underline{u}}}
|\underline{L}^{2}\Omega\phi|^{2}=\int_{C_{u_{0}}}|\slashed{\nabla}\underline{L}\Omega\phi|^{2}+
\iint_{\mathcal{D}}\phi\underline{L}\Omega\phi\underline{L}^{2}\Omega\phi+\iint_{\mathcal{D}}\underline{L}
\phi\Omega\phi\underline{L}^{2}\Omega\phi\\
+\iint_{\mathcal{D}}\frac{1}{r^{2}}(\underline{L}\Omega\phi-L\Omega\phi)\underline{L}^{2}\Omega\phi+
\iint_{\mathcal{D}}\frac{1}{r^{3}}
\slashed{\Delta}\Omega\phi \underline{L}^{2}\Omega\phi\\
=\int_{C_{u_{0}}}|\slashed{\nabla}\underline{L}\Omega\phi|^{2}+T_{1}+T_{2}+T_{3}+T_{4}
\end{align*}
By Proposition 4.1 and bootstrap assumption,
\begin{align*}
 T_{1}\lesssim\|\phi\|_{L^{\infty}(S_{\underline{u},u})}\int_{0}^{\underline{u}}\|\underline{L}\Omega\phi\|_
{L^{2}(\underline{C}_{\underline{u}^{\prime}})}\|\underline{L}^{2}\Omega\phi\|
_{L^{2}(\underline{C}_{\underline{u}^{\prime}})}d\underline{u}^{\prime}\\
\lesssim(\delta^{1/4}M)\delta(\delta^{1/2}M)(M)\lesssim\delta^{7/4}M^{3}
\end{align*}
\begin{align*}
 T_{2}\lesssim\int_{0}^{\underline{u}}\|\underline{L}\phi\|_{L^{4}(\underline{C}_{\underline{u}^{\prime}})}
\|\Omega\phi\|_{L^{4}(\underline{C}_{\underline{u}^{\prime}})}\|\underline{L}^{2}\Omega\phi\|
_{L^{2}(\underline{C}_{\underline{u}^{\prime}})}
d\underline{u}^{\prime}\lesssim\delta(\delta^{1/4}M)(\delta^{1/4}M)M\lesssim\delta^{3/2}M^{3}
\end{align*}
also
\begin{align*}
 T_{3}\lesssim\int_{0}^{\underline{u}}\|\underline{L}\Omega\phi\|_{L^{2}(\underline{C}_{\underline{u}^{\prime}})}
\|\underline{L}^{2}\Omega\phi\|_{L^{2}(\underline{C}_{\underline{u}^{\prime}})}d\underline{u}^{\prime}
+(\int_{u_{0}}^{u}\|L\Omega\phi\|^{2}_{L^{2}(C_{u^{\prime}})}du^{\prime})^{1/2}(\int_{0}^{\underline{u}}
\|\underline{L}^{2}\Omega\phi\|^{2}
_{L^{2}(\underline{C}_{\underline{u}^{\prime}})}d\underline{u}^{\prime})^{1/2}\\
\lesssim\delta^{3/2}M^{2}+\delta^{1/2}M^{2}\lesssim\delta^{1/2}M^{2}
\end{align*}
Similarly,
\begin{align*}
 T_{4}\lesssim\int_{0}^{\underline{u}}\|\slashed{\Delta}\Omega\phi\|_{L^{2}(\underline{C}_{\underline{u}^{\prime}})}
\|\underline{L}^{2}\Omega\phi\|_{L^{2}(\underline{C}_{\underline{u}^{\prime}})}d\underline{u}^{\prime}\lesssim\delta 
M^{2}
\end{align*}
So we obtain:
\begin{align*}
 \int_{\underline{C}_{\underline{u}}}|\underline{L}^{2}\Omega\phi|^{2}\lesssim I^{2}
+\delta^{1/2}M^{3}
\end{align*}
i.e.
\begin{align}
 \underline{F}_{3}(u,\underline{u})\lesssim I+\delta^{1/4}M^{3/2}
\end{align}
\subsection{End of the Bootstrap Argument} Combining (39), (45), (46) and (47) we obtain:
\begin{align}
 M\lesssim I+\delta^{1/4}M^{3/2}
\end{align}
Choosing $\delta$ sufficiently small depending on the quantity $I$ together with the bootastrap assumption:
\begin{align*}
 M\leq AI
\end{align*}
 we obtain:
\begin{align}
 M\leq \frac{A}{2}I
\end{align}
Therefore $\mathcal{A}$ 
is both an open and closed subset of $\mathcal{U}$, then is $\mathcal{U}$ itself.

This completes the proof of $\textbf{Main A priori Estimates}$.
\subsection{Higher Order Estimates} For the estimates of higher order derivatives, we just use the induction to prove
 it, since we have established the estimates for the lower derivatives up to the 3rd order. With the definitions:
\begin{align*}
 E_{k}(u,\underline{u})=\|L\Omega^{k-1}\phi\|_{L^{2}(C_{u})}+\delta^{-1/2}\|\Omega^{k}\phi\|_{L^{2}(C_{u})}\\
\underline{E}_{k}(u,\underline{u})=\|\Omega^{k}\phi\|_{L^{2}(\underline{C}_{\underline{u}})}
+\delta^{-1/2}\|\underline{L}\Omega^{k-1}\phi\|_{L^{2}(\underline{C}_{\underline{u}})}
\end{align*}
and
\begin{align*}
 F_{k}(u,\underline{u})=\delta\|L^{2}\Omega^{k-2}\phi\|_{L^{2}(C_{u})}\quad
\underline{F}_{k}(u,\underline{u})=\|\underline{L}^{2}\Omega^{k-2}\phi\|_{L^{2}(\underline{C}_{\underline{u}})}
\end{align*}
Then if the initial data satisfy 
\begin{align*}
 \sum_{i=1}^{n+2}E_{i}(u_{0},\delta)+\sum_{j=2}^{n+2}F_{j}(u_{0},\delta)\leq I
\end{align*}
we have:
\begin{align*}
 E_{n+2}(u,\underline{u})+\underline{E}_{n+2}(u,\underline{u})+F_{n+2}(u,\underline{u})+\underline{F}_{n+2}
(u,\underline{u})
\lesssim_{n} I
\end{align*}
Similar as Proposition 4.1,
\begin{align*}
  \delta^{-1/4}\|\Omega^{n-1}\phi\|_{L^{\infty}(S_{\underline{u},u})}
+\delta^{1/2}\|L\Omega^{n}\phi\|_{L^{4}(S_{\underline{u},u})}+\delta^{-1/4}\|\Omega^{n+1}\phi\|
_{L^{4}(S_{\underline{u},u})}
+\delta^{-1/4}\|\underline{L}\Omega^{n}\phi\|_{L^{4}(S_{\underline{u},u})}\\
+\delta^{1/2}\|L\Omega^{n-1}\phi\|_{L^{4}(S_{\underline{u},u})}+\delta^{-1/4}\|\Omega^{n}\phi\|
_{L^{4}(S_{\underline{u},u})}
+\delta^{-1/4}\|\underline{L}\Omega^{n-1}\phi\|_{L^{4}(S_{\underline{u},u})}\lesssim_{n} I
\end{align*}
\section{Existence of Solutions}  By the a priori estimates, we can show that (1) with data prescribed on $C_{u_{0}}$ 
where $u_{0}\leq\underline{u}\leq\delta$ can be solved all the way up to $t=-1$.

    We use the local existence result of Alan. D. Rendall \cite{Ren}, which states that there exists a solution around 
$S_{0,u_{0}}$, say, defined in the region enclosed by $\underline{C}_{0}$, $C_{u_{0}}$ and $t=u_{0}+\epsilon$ with 
$\epsilon\ll\delta$. Thanks to the a priori Estimates, the solution and its derivatives are bounded
on $t=u_{0}+\epsilon$ by the initial data. Therefore, we can solve a Cauchy problem with data prescribed on 
$t=u_{0}+\epsilon$ to construct a solution in the future domain of dependence of $t=u_{0}+\epsilon$ whose boundary 
contains two null hypersurfaces $C_{u_{0}+\epsilon}$ and $\underline{C}_{\epsilon}$. Now we have two
characteristic problem: for the first one, the data is prescribed on $\underline{C}_{0}$ and $C_{u_{0}+\epsilon}$; 
for the second one, the data is prescribed on $C_{u_{0}}$ and $\underline{C}_{\epsilon}$. We can use Rendall's local 
existence result again to solve them around $S_{0,u_{0}+\epsilon}$ and $S_{\epsilon,u_{0}}$.
In this way, we can actually push the solution to $t=u_{0}+\epsilon+\epsilon^{\prime}$ with another small 
$\epsilon^{\prime}$. Then we can repeat the above process to push the solution all the way to $t=u_{0}+\delta$, 
and then to $t=-1$. Actually, from the second step, since we the main A priori estimate, the the length of the 
interval where the solution exists is the same. Therefore we can finally push the solution to $t=-1$.

\section{Construction of Cauchy Data, Final Conclusions}

$\textbf{Proposition 6.1}$ \emph{Assume we have bound on $E_{i}(u_{0})$ with $i=1,2,...,n+2$ and $F_{j}(u_{0})$ with 
$j=2,3,...,n+2$ for some fixed $n\geq 10$. Then for $k=4,...,n$, we have}
\begin{align*}
 \|\Omega^{k-4}\phi\|_{L^{\infty}(\underline{C}_{\delta}}+...+\|\Omega^{k-1}\phi\|_{L^{\infty}(\underline{C}_{\delta})}+
\|\underline{L}\Omega^{k-3}\phi\|_{L^{\infty}(\underline{C}_{\delta})}+
\|\underline{L}^{2}\Omega^{k-4}\phi\|_{L^{\infty}(\underline{C}_{\delta})}\lesssim\delta^{1/4}\\
\|L\Omega^{k-3}\phi\|_{L^{\infty}(\underline{C}_{\delta})}+\|L^{2}\Omega^{k-4}\phi\|_{L^{\infty}
(\underline{C}_{\delta})}\lesssim\delta^{1/4}
\end{align*}
$Proof$. The estimate for $\|\Omega^{k-1}\phi\|_{L^{\infty}(\underline{C}_{\delta})}$ comes from the property 
 before section 5. For the $L$ and $\underline{L}$ derivatives, we just consider $L\phi$ and $\underline{L}\phi$, 
the higher order derivatives are similar.
To start, we write equation (15) in two different forms:
\begin{align}
 L(\underline{L}\phi)=a(\underline{L}\phi)+\bar{b}+\frac{1}{r}(L\phi)
\end{align}
and also
\begin{align}
 \underline{L}(L\phi)=-a(L\phi)+\bar{b}-\frac{1}{r}(\underline{L}\phi)
\end{align}
where
\begin{align*}
 \bar{b}=\frac{1}{r^{2}}\slashed{\Delta}\phi\mp|\phi|^{k-1}\phi
\end{align*}
By Gronwall's inequality, we have, since $L\phi$ vanishes near $S_{\delta,u_{0}}$:
\begin{align*}
 \|L\phi\|_{L^{\infty}(S_{\underline{u},u})}\lesssim\int_{u_{0}}^{u}(
\|\underline{L}\phi\|_{L^{\infty}(S_{\underline{u},u^{\prime}})}+\|b\|_{L^{\infty}(S_{\underline{u},u^{\prime}})})
du^{\prime}
\end{align*}
(Here $\underline{u}$ is very close to $\delta$) and also
\begin{align*}
 \|\underline{L}\phi\|_{L^{\infty}(S_{\underline{u},u})}\lesssim
\int_{0}^{\underline{u}}(\|L\phi\|_{L^{\infty}(S_{\underline{u}^{\prime},u})}+\|b\|
_{L^{\infty}(S_{\underline{u}^{\prime},u})})d\underline{u}^{\prime}
\end{align*}
Defining
\begin{align*}
 A(u,\underline{u})=\sup_{\underline{u}^{\prime}\in[0,\underline{u}]}\|L\phi\|_{L^{\infty}(S_{\underline{u}^{\prime},u})}
\\B(u,\underline{u})=\sup_{u^{\prime}\in[u_{0},u]}\|\underline{L}\phi\|_{L^{\infty}(S_{\underline{u},u^{\prime}})}
\end{align*}
Since
\begin{align*}
 \|b\|_{L^{\infty}(S_{\underline{u},u})}\lesssim\delta^{1/4}
\end{align*}
we obtain:
\begin{align}
 A(u,\underline{u})\lesssim B(u,\underline{u})+\delta^{1/4}\\
B(u,\underline{u})\lesssim\delta A(u,\underline{u})+\delta^{5/4}
\end{align}
where we have used the fact that
\begin{align*}
 A(u,\underline{u}^{\prime})\lesssim A(u,\underline{u})\\
B(u^{\prime},\underline{u})\lesssim B(u,\underline{u})
\end{align*}
Substituting (53) in (52),
\begin{align*}
 A(u,\underline{u})\lesssim \delta A(u,\underline{u})+\delta^{1/4}
\end{align*}
So if we choose $\delta$ sufficiently small, we obtain:
\begin{align*}
 A(u,\underline{u})\lesssim\delta^{1/4}
\end{align*}
Back to (53),
\begin{align*}
 B(u,\underline{u})\lesssim\delta^{5/4}
\end{align*}
Then the proposition follows. $\qed$

So we obtain from the above proposition that the data on $\underline{C}_{\delta}$ induced from the solution 
are small in energy norms. 
Note also that we lose one derivative when we integrate the propagation equation because of the term 
$\slashed{\Delta}\phi$.

Now we can construct our Cauchy data, which is similar to \cite{WY1}, and whose energy is larger than $E_{0}$

We now choose a Cauchy hypersurface $\Sigma=\{t=u_{0}+\delta\}$. Let \\
$\Sigma_{1}=\Sigma\bigcap(\textrm{Region 1}\bigcup\textrm{Region 2})$
and $\Sigma_{2}=\Sigma-\Sigma_{1}$. By the above proposition, we know that there is an an annular region $E$ bounded 
by $S_{\delta, u_{0}}$ and a smaller sphere near $S_{\delta,u_{0}}$ in $\Sigma_{1}$, on which the solution is small. 
Then by a Whitney extension theorem established by Fefferman \cite{Fef} and using a cut off function, we can extend the 
data $(\phi^{(0)}_{1},\phi^{(1)}_{1})$ on $\Sigma_{1}$ to the whole $\Sigma$ with the following properties:
\begin{align*}
 (\phi^{(0)},\phi^{(1)})|_{\Sigma_{1}}=(\phi^{(0)}_{1},\phi^{(1)}_{1})\\
(\phi^{(0)},\phi^{(1)})|_{\{x\in\Sigma_{2}|\textrm{dist}(x,\Sigma_{1})\geq 1\}}=(0,0)\\
\|\phi^{(0)}\|_{L^{\infty}(\{x\in\Sigma_{2}|\textrm{dist}(x,\Sigma_{1})\leq 1\})},\\
\|(\partial^{k-1}\phi^{(0)},\partial^{k-2}\phi^{(1)})\|_{L^{\infty}(\{x\in\Sigma_{2}|\textrm{dist}(x,\Sigma_{1})
\leq 1\})}\lesssim\delta^{1/4}
\quad \textrm{for}\quad k=2,3,...,n
\end{align*}
where we denote by $\partial^{k-1}\phi_{0}$ and $\partial^{k-2}\phi_{1}$ the derivatives appearing in Proposition 6.1.

Therefore, according to the small data theory, we obtain a solution $\phi$ in Region 4. See \cite{Joh}. In particular, 
the energy flux on $C^{+}_{u_{0}}$ induced from the solution in Region 4 are small.

     We now have the data on $\underline{C}_{\delta}$ and $C^{+}_{u_{0}}$. They are past boundaries of Region 3. We can
then solve this small data problem in Region 3. Together with the solutions constructed in other regions, this 
completes the construction of the whole solution.

     Next, we must show that the energy of the initial data is larger than $E_{0}$, provided that $\delta$ is suitably 
small. Recall the definition:
\begin{align*}
 E(\phi(t))=\frac{1}{2}\int_{\Sigma_{t}}|\partial_{t}\phi|^{2}+|\nabla_{x}\phi|^{2}\pm\frac{2}{k+1}|\phi|^{k+1}dx
\end{align*}
By Proposition 4.1 we have an $L^{\infty}$ bound for $\phi$, and note also that the solution on $\Sigma_{u_{0}+\delta}$
is compactly supported in an annular domain of size $\delta$. So the potential energy will be very small. We must prove 
that the kinetic energy is large.

    To do this, we use the energy identity on the domain bounded by $\underline{C}_{0}$, $C_{u_{0}}$ and 
$\{t=u_{0}+\delta\}$. Since by (9), the spacetime integral is clearly small (depending on $\delta$), and the solution
 vanishes on $\underline{C}_{0}$, so the energy considered is comparable to the energy on $C_{u_{0}}$, which can be 
larger than $E_{0}$, if we choose $\psi_{0}$ properly. This completes the proof of the main conclusion.

\section{2-D Case--A Sketch}
Our method can also be used to deal with the equation in $\mathbb{R}\times\mathbb{R}^{2}$. We will use the equation:
\begin{align}
 \Box\phi=-\phi e^{\phi^{2}}\quad \textrm{in}\quad \mathbb{R}^{2}\times\mathbb{R}
\end{align}
where
\begin{align*}
 \Box=-\partial^{2}_{tt}+\Delta_{x}
\end{align*}
as an example, the general case can be dealt similarly. The energy associated to (54) is
\begin{align}
 E(\phi)=\int_{\mathbb{R}^{3}}(\frac{1}{2}|\partial_{t}\phi(t,x)|^{2}+\frac{1}{2}|\nabla_{x}\phi(t,x)|^{2}-
e^{\phi(t,x)^{2}}+1)dx
\end{align}

The equation (54) is the focusing, energy super-critical nonlinear wave equation, see \cite{IMM} and \cite{Str2} for 
a reference. There are few results about this equation in the focusing case.\\

Now in $\mathbb{R}^{2+1}$, we have the rectangular coordinates $(t,x_{1},x_{2})$ as well as null-polar coordinates 
$(u,\underline{u},\theta)$. The geometric setting and the energy identity is almost the same as in the case 
$\mathbb{R}^{3+1}$. Now the rotation vectorfield is
\begin{align*}
 \Omega=\partial_{\theta}=x_{1}\partial_{2}-x_{2}\partial_{1}
\end{align*}
 $S_{\underline{u},u}$ is a circle, and the energy currents associated to the deformation tensors are
\begin{align}
 K^{\Omega}=0\quad K^{L}=\frac{1}{2r}(|\slashed{\nabla}\phi|^{2}+L\phi\underline{L}\phi)\quad
K^{\underline{L}}=-\frac{1}{2r}(|\slashed{\nabla}\phi|^{2}+L\phi\underline{L}\phi)
\end{align}

Since now we in $\mathbb{R}^{2}\times\mathbb{R}$, the Sobolev inequalities will be different. Actually, we have:

$\textbf{Lemma 7.1}$ \emph{For a smooth function $f$ on the circle $S_{\underline{u},u}$, we have}
\begin{align*}
 \sup_{S_{\underline{u},u}}|f|\leq|u|^{1/2}(\int_{S_{\underline{u},u}}|\slashed{\nabla}f|^{2}d\mu_{\slashed{g}})^{1/2}
+|u|^{-1/2}(\int_{S_{\underline{u},u}}|f|^{2}d\mu_{\slashed{g}})^{1/2}\\
\int_{S_{\underline{u},u}}|f|^{6}d\mu_{\slashed{g}}\leq(\int_{S_{\underline{u},u}}|f|^{4}d\mu_{\slashed{g}})
\{|u|\int_{S_{\underline{u},u}}|\slashed{\nabla}f|^{2}d\mu_{\slashed{g}}+|u|^{-1}\int_{S_{\underline{u},u}}|f|^{2}
d\mu_{\slashed{g}}\}
\end{align*}
This lemma can be proved by using the isoperimetric inequality on circles.

Also we have:

$\textbf{Lemma 7.2}$ \emph{Let $\phi$ be a smooth function on $C_{u}$ vanishing on $S_{0,u}$, then we have}
\begin{align*}
 |u|^{1/4}\|\phi\|_{L^{4}(S_{\underline{u},u})}\lesssim\|L\phi\|^{1/2}_{L^{2}(C_{u})}
(\|\phi\|^{1/2}_{L^{2}(C_{u})}+|u|^{1/2}\|\slashed{\nabla}\phi\|^{1/2}_{L^{2}(C_{u})})\\
\|\phi\|_{L^{2}(S_{\underline{u},u})}\lesssim\|L\phi\|^{1/2}_{L^{2}(C_{u})}\|\phi\|^{1/2}_{L^{2}(C_{u})}
\end{align*}
 $\textbf{Lemma 7.3}$ \emph{Let $\phi$ be a smooth function on $\underline{C}_{\underline{u}}$, we have the following 
estimates}:
\begin{align*}
 |u|^{1/4}\|\phi\|_{L^{4}(S_{\underline{u},u})}\lesssim|u_{0}|^{1/4}\|\phi\|_{L^{4}(S_{\underline{u},u_{0}})}
+\|\underline{L}\phi\|^{1/2}_{L^{2}(\underline{C}_{\underline{u}})}(\|\phi\|^{1/2}_{L^{2}(\underline{C}_{\underline{u}})}
+\||u^{\prime}|\slashed{\nabla}\phi\|^{1/2}_{L^{2}(\underline{C}_{\underline{u}})})\\
\|\phi\|_{L^{2}(S_{\underline{u},u})}\lesssim\|\phi\|_{L^{2}(S_{\underline{u},u_{0}})}+
\|\underline{L}\phi\|^{1/2}_{L^{2}(\underline{C}_{\underline{u}})}\|\phi\|^{1/2}_{L^{2}(\underline{C}_{\underline{u}})}
\end{align*}
The proof of Lemma 7.2 and Lemma 7.3 can be found in \cite{WY1}.

We construct the characteristic initial data in the same way as 3-D case. The main difference in 2-D case is that we 
only need the first and the second derivatives of the solution to close the bootstrap, this is because the Sobolev 
inequalities involve one less derivative in 2-D case. Namely, we 
define the following norms which has one less derivative than \cite{WY1}:
\begin{align}
 E_{1}(u,\underline{u})=\|L\phi\|_{L^{2}(C_{u})}+\delta^{-\frac{1}{2}}\|\Omega\phi\|_{L^{2}(C_{u})}\\\notag
\underline{E}_{1}(u,\underline{u})=\|\Omega\phi\|_{L^{2}(\underline{C}_{\underline{u}})}+
\delta^{-\frac{1}{2}}\|\underline{L}\phi\|_{L^{2}(\underline{C}_{\underline{u}})},\\\notag
E_{2}(u,\underline{u})=\|L\Omega\phi\|_{L^{2}(C_{u})}+\delta^{-\frac{1}{2}}\|\Omega^{2}\phi\|_{L^{2}(C_{u})},\\\notag
\underline{E}_{2}(u,\underline{u})=\|\Omega^{2}\phi\|_{L^{2}(\underline{C}_{\underline{u}})}+\delta^{-\frac{1}{2}}
\|\underline{L}\Omega\phi\|
_{L^{2}(\underline{C}_{\underline{u}})},\\\notag
\end{align}
and also
\begin{align}
 F_{2}(u,\underline{u})=\delta\|L^{2}\phi\|_{L^{2}(C_{u})},\\\notag
\underline{F}_{2}(u,\underline{u})=\|\underline{L}^{2}\phi\|_{L^{2}(\underline{C}_{\underline{u}})}
\end{align}
then we obtain the following result which is similar to that of 3-D case:

$\textbf{Theorem 7.1}$ \emph{If $\delta$ is sufficiently small, for all characteristic initial data of (54) and all positive 
real number $I$ which satisfy}
\begin{align}
E_{1}(u_{0},\delta)+E_{2}(u_{0},\delta)+F_{2}(u_{0},\delta)\leq I
\end{align}
\emph{there is a constant $C(I)$ depending only on $I$, so that}
\begin{align}
 \sum_{i=1}^{2}[E_{i}(u,\underline{u})+\underline{E}_{i}(u,\underline{u})]+F_{2}(u,\underline{u})+\underline{F}
_{2}(u,\underline{u})\leq C(I)
\end{align}
Once we establish the above result, the following steps are exactly the same as 3-D case.

\section*{Acknowledgements}
\thispagestyle{empty}

The author would like to thank Jonas L\"uhrmann for his helpful comments and suggestions on introduction and the
construction of Cauchy data; Prof. Demetrios Christodoulou and Prof. Pin Yu for communications on $short$ $pulse$
method; and Prof. Ping Zhang for his long standing encouragement. 
This work is supported by the ERC Advanced Grant No. 246574 ``Partial Differential Equations of Classical Physics'' 
directed by Prof. Christodoulou.
\clearpage

\bibliographystyle{amsplain}
\bibliography{myreference}\vspace{12mm}

\noindent Department of Mathematics,\\ 
ETH Zurich,\\ 
R\"amistrasse 101, 8092 Zurich,\\
Switzerland\\
\emph{Email}: shuang.miao@math.ethz.ch\\[.2cm]
and\\[.2cm]
Academy of Mathematics and Systems Sciences,\\
Chinese Academy of Science,\\ 
Zhongguancun East Road 55, 100190 Beijing,\\
China\\
\emph{Email}: miaoshuang@amss.ac.cn\\

\end{document}